\documentclass{article}
\usepackage{amsfonts}

\input{tcilatex}
\input{tcilatex}
\begin{document}

\begin{center}
$\mathbf{p}$\textbf{-Laplacian Fractional Sturm-Liouville Problem for
Diffusion Operator via Impulsive Condition}

\bigskip

$^{1}$Funda METIN TURK, $^{2}$Erdal BAS

\bigskip $^{1}$\textit{Department of Mathematics, Faculty of Science, Firat
University, Elazig, 23119, Turkey}

$^{2}$\textit{Department of Mathematics, Faculty of Science, Firat
University, Elazig, 23119, Turkey}

\textit{e-mail}$^{1}$\textit{: fnd-44@hotmail.com}

\textit{e-mail}$^{2}$\textit{: erdalmat@yahoo.com}

\bigskip
\end{center}

\textbf{Abstract. }In this study, the existence results of solution is given
for fractional $p$-Laplacian Stum-Liouville problem for diffusion operator
of order with impulsive conditions. The derivatives are described in
Riemann-Liouville and Caputo sense. The Riemann-Liouville integral operator
is used to acquire the integral representation of solution. The existence of
solution is demonstrate via Schaefer fixed point theorem.



\bigskip

\textbf{1. Introduction}

\bigskip

The last half of the past century has witnessed to both intensive
improvement of the theory of diffrential equations involving derivatives
fractional order and the applications such as physics, control systems,
polymer rheology, aerodynamics and other areas. Fractional differential
equations have been constantly drawing interest of many autors. The
attention in the study of fractional differential equations is based upon
the fact that fractional calculus service as an great tool in common usage
for the applications of such constructions in various sciences, the
description of properties of diverse materials, processes and important part
of the physical mathematics and also a large part of the literature is
related to fractional differential equations. In consequence, the fractional
order models is more factual and useful than the integer order models. For
more information and applications about fractional calculus, see
[1-10,14-21] and references therein. At the same time, the fractional $p$%
-Laplacian operator appear naturally in the applied sciences and is
extensively used in the mathematical modeling of physical and natural
phenomena, blood flow problems, turbulent filtration in porous media,
rheology, modeling of viscoplasticity, mechanics, material science and many
other related fields. Therefore, a continuous increasing attention has been
shown towards problems involving the fractional $p$-Laplacian operator and
on the existence of solutions for this problem. But there is no known study
about the existence of solutions for fractional $p$-Laplacian
Sturm-Liouville problem thus this paper is a main study on literature. For
examples and details, see [22-26]

Impulsive differential equations have arisen as a significant area for
applied sciences in recent years. Impulsive differential equations are
accepted as significant mathematical devices to make many real world
problems plausible in applied sciences. There is a great deal of study for
boundary value problems of impulsive differential equations of integer order
in the literature. On the other hand, there is very little known about
impulsive boundary value problems for fractional order and many aspect of
these problems are yet to be discovered.

Recently, there has been too much attention on the existence of solutions
for impulsive boundary-value problems for fractional differential equations
by means of techniques (fixed point theorems, Banach contraction mapping
principle, etc.). This subject has been studied in the various papers
[27--29,32-37]. For example, in [32], Yuansheng Tian and Zhanbing Bai
discussed the existence results for the three-point impulsive boundary value
problem involving fractional differential equations given by%
\[
^{C}D^{\alpha }u\left( t\right) =f\left( t,u\left( t\right) \right) ,\text{
\ }0<t<1,\text{ }t\neq t_{k},\text{ }k=1,2,...,p, 
\]%
\[
\left. \Delta u\right\vert _{t=t_{k}}=I_{k}\left( u\left( t_{k}\right)
\right) ,\text{ }\left. \Delta u^{\prime }\right\vert _{t=t_{k}}=\bar{I}%
_{k}\left( u\left( t_{k}\right) \right) ,\text{ }k=1,2,..,p, 
\]%
\[
u\left( 0\right) +u^{\prime }\left( 0\right) =0,\text{ }u\left( 1\right)
+u^{\prime }\left( \xi \right) =0, 
\]%
where $^{C}D^{a}$ is the Caputo fractional derivative, $q\in R,$ $1<q\leq 2,$
$f:\left[ 0,1\right] \times R\rightarrow R$ is a continuous function, $I_{k},%
\bar{I}_{k}:R\longrightarrow R,$ $\xi \in \left( 0,1\right) ,\xi \neq
t_{k},k=1,2,...,p$ and $\left. \Delta u\right\vert _{t=t_{k}}=u\left(
t_{k}^{+}\right) -u\left( t_{k}^{-}\right) ,\left. \Delta u^{\prime
}\right\vert _{t=t_{k}}=u^{\prime }\left( t_{k}^{+}\right) -u^{\prime
}\left( t_{k}^{-}\right) ,$ $u\left( t_{k}^{+}\right) $ and $u\left(
t_{k}^{-}\right) $ shows the right and left-hand limit of the function $%
u\left( t\right) $ at $t=t_{k}$, and the sequences $\left\{ t_{k}\right\} $
satisfy that $0=t_{0}<t_{1}<...<t_{p}<t_{p+1}=1,p\in 
\mathbb{N}
.$

Ravi P. Agarwal, Mouffak Benchohra and Boualem Attou Slimani [34]
investigated the existence and uniqueness of solutions for the initial value
problems for fractional order differential equations as the following form%
\begin{eqnarray*}
^{C}D^{\alpha }y\left( t\right) &=&f\left( t,y\right) ,\text{ \ for each }%
t\in J=\left[ 0,T\right] , \\
t &=&t_{k},\text{ }k=1,2,..,m\text{, }1<\alpha \leq 2
\end{eqnarray*}%
\[
\left. \Delta y\right\vert _{t=t_{k}}=I_{k}\left( y\left( t_{k}^{-}\right)
\right) ,\text{ }\left. \Delta y^{\prime }\right\vert _{t=t_{k}}=\bar{I}%
_{k}\left( y\left( t_{k}^{-}\right) \right) ,\text{ }k=1,2,..,m, 
\]%
\[
y\left( 0\right) =y_{0}\left( 0\right) ,\text{ }y^{\prime }\left( 0\right)
=y_{1}, 
\]%
where $^{C}D^{a}$ is the Caputo fractional derivative, $f:J\times
R\rightarrow R$ is a continuous function, $I_{k},\bar{I}_{k}:R%
\longrightarrow R,$ $k=1,2,...,m$ and $y_{0},y_{1}\in R,$ $%
0=t_{0}<t_{1}<...<t_{m}<t_{m+1}=T,$ $\left. \Delta y\right\vert
_{t=t_{k}}=y\left( t_{k}^{+}\right) -y\left( t_{k}^{-}\right) ,$ $y\left(
t_{k}^{+}\right) =\lim\limits_{h\rightarrow 0^{+}}y\left( t_{k}+h\right) $
and $y\left( t_{k}^{-}\right) =\lim\limits_{h\rightarrow 0^{-}}y\left(
t_{k}+h\right) $ at $t=t_{k},$ $k=1,2,...,m$.

Moreover, Jie Zhou and Meiqiang Feng [37] study fractional Sturm-Liouville
problem with impulsive condition.

The object of this study is to develop main parts of Sturm-Liouville theory
for the $p$-Laplacian and is to continue this study by giving several
existence results for fractional $p$-Laplacian Sturm-Liouville problem
having diffusion operator with impulsive conditions and is to prosper the
theoretical knowledge of the above. Therefore we analyze the following
fractional $p$-Laplacian Sturm-Liouville problem having diffusion operator
with impulsive conditions%
\begin{equation}
-D_{0,+}^{\beta }\phi _{p}\;^{C}D_{0,+}^{\alpha }y\left( t\right) +\left(
2\lambda p\left( t\right) +q\left( t\right) \right) y\left( t\right) =0, 
\tag{1}
\end{equation}%
\begin{equation}
\Delta y\mid _{t=t_{k}}=I_{k}\left( y\left( t_{k}\right) \right) ,\text{ }%
\Delta y^{\prime }\mid _{t=t_{k}}=I_{k}^{\ast }\left( y\left( t_{k}\right)
\right) ,\text{ }t_{k}\in \left( 0,\pi \right) ,\text{ }k=1,2,..,n,  \tag{2}
\end{equation}%
\begin{equation}
y\left( 0\right) +y^{\prime }\left( 0\right) =0,\text{ \ }y\left( \pi
\right) +y^{\prime }\left( \pi \right) =0,  \tag{3}
\end{equation}

\noindent where $D_{0+}^{\beta }$ Riemann-Liouville, $^{C}D_{0+}^{a}$ is the
Caputo fractional derivative, $p\in W,q\in L^{2}\left[ 0,\pi \right] $, $%
I_{k},I_{k}^{\ast }:%
\mathbb{R}
\longrightarrow 
\mathbb{R}
$ are real-valued continuous functions, $\left. \Delta y\right\vert
_{t=t_{k}}=y\left( t_{k}^{+}\right) -y\left( t_{k}^{-}\right) $, $\ y\left(
t_{k}^{+}\right) =\lim\limits_{h\longrightarrow 0^{+}}y\left( t_{k}+h\right) 
$, $y\left( t_{k}^{-}\right) =\lim\limits_{h\longrightarrow 0^{-}}y\left(
t_{k}+h\right) $ at $t=t_{k}$, $k=1,2,..,n$, $%
0=t_{0}<t_{1}<...<t_{n}<t_{n+1}=\pi $, $\Delta y^{\prime }\mid _{t=t_{k}}$
has a similar meaning for $y^{\prime }\left( t\right) .$ $\phi _{p}\left(
s\right) =\left\vert s\right\vert ^{p-2}s,p>1$. Obviously, $\phi _{p}$ is
invertible and its inverse operator is $\phi _{q}$, where $q>1$ is a
constant such that $\frac{1}{p}+\frac{1}{q}=1$.

Some necessary notations, definitions and lemmas are given in Section 2 . We
establish a theorem on existence of solution for (1)-(3) problem by using
Schaefer's fixed point theorem for $p$-Laplacian fractional Sturm-Liouville
problem via diffusion operator in Section 3.

\bigskip

\textbf{2.} \textbf{Preliminaries}

We give some material related to fractional calculus theory. For more
details about this field, see [4,5,8].

Considering the following space%
\begin{eqnarray*}
PC\left( J,R\right) &=&\left\{ y:J\longrightarrow R:y\in C\left( \left(
t_{k},t_{k+1}\right] ,R\right) ,k=0,...,n+1\text{ and }\right. \\
&&\left. \text{there exist }y\left( t_{k}^{-}\right) \text{ and }y\left(
t_{k}^{+}\right) \text{, }k=1,2,...,n,y\left( t_{k}^{-}\right) =y\left(
t_{k}\right) \right\}
\end{eqnarray*}%
where $J=\left[ 0,\pi \right] $. $PC\left( J,R\right) $ is a Banach space
with the norm%
\[
\left\Vert y\right\Vert _{PC}=\sup_{t\in J}\left\vert y\left( t\right)
\right\vert . 
\]

\bigskip

\textbf{Definition 1.} The left and right-sided Riemann-Liouville integrals
of order $\alpha $ are given by [4].

\bigskip

\textbf{Definition 2. }[4]The left and right-sided Riemann-Liouville
derivatives are defined as respectively, $0<\alpha <1,$%
\[
\left( D_{a+}^{\alpha }f\right) \left( r\right) =D\left( I_{a+}^{1-\alpha
}f\right) \left( r\right) \text{ \ \ \ }r>a, 
\]%
\[
\left( D_{b-}^{\alpha }f\right) \left( r\right) =-D\left( I_{b-}^{1-\alpha
}f\right) \left( r\right) \text{ \ \ }r<b. 
\]%
Similar formulas give the left and right-sided Caputo derivatives of order $%
\alpha $:%
\[
\left( ^{C}D_{a+}^{\alpha }f\right) \left( r\right) =\left( I_{a+}^{1-\alpha
}Df\right) \left( r\right) \text{ \ \ \ }r>a 
\]%
\[
\left( ^{C}D_{b-}^{\alpha }f\right) \left( r\right) =\left( I_{b-}^{1-\alpha
}\left( -D\right) f\right) \left( r\right) \text{ \ \ }r<b 
\]

\bigskip

\textbf{Definition 3. }[30] If $K$ is a compact metric space then a subset $%
F\subset C\left( K\right) $ of the space of continuous functions on $K$
equipped with the uniform distance, is compact if and only if it is closed,
bounded and equicontinuous.

\bigskip

\textbf{Lemma 4. }[31] Let $\alpha >0$. Then the differential equation%
\[
^{C}D^{\alpha }h\left( t\right) =0, 
\]%
has solution $h\left( t\right) =c_{0}+c_{1}t+c_{2}t^{2}+...+c_{n}t^{n-1}$, $%
c_{i}\in 
\mathbb{R}
,$ $i=0,1,2,...,n,$ $n=\left[ \alpha \right] +1$.

\bigskip

\textbf{Lemma 5. }[31] Let $\alpha >0$. Then%
\[
I^{\alpha }{}^{C}D^{\alpha }h\left( t\right) =h\left( t\right)
+c_{0}+c_{1}t+c_{2}t^{2}+...+c_{n}t^{n-1}, 
\]%
for some $c_{i}\in 
\mathbb{R}
,$ $i=0,1,2,...,n,$ $n=\left[ \alpha \right] +1$.

\bigskip

\textbf{Lemma 6. }[4] $\func{Re}\left( \alpha \right) >0,$ $n=\func{Re}%
\left( \alpha \right) +1$ and let $f_{n-\alpha }\left( x\right) =\left(
I_{b-}^{n-\alpha }f\right) \left( x\right) $ be the fractional integral of
order $n-\alpha .$

\textbf{a)} If $1\leq p\leq \infty $ and $f\left( x\right) \in
I_{b-}^{\alpha }\left( L_{p}\right) $, then%
\[
\left( I_{b-}^{\alpha }D_{b-}^{\alpha }f\right) \left( x\right) =f\left(
x\right) 
\]%
where%
\[
I_{b-}^{\alpha }\left( L_{p}\right) =\left \{ f:f=I_{b-}^{\alpha }\varphi ,%
\text{ }\varphi \in L_{p}\left( a,b\right) \right \} . 
\]

\textbf{b)} If $f\left( x\right) \in L_{1}\left( a,b\right) $ and $%
f_{n,...,\alpha }\left( x\right) \in AC^{n}\left[ a,b\right] ,$ then the
equality%
\[
\left( I_{b-}^{\alpha }D_{b-}^{\alpha }f\right) \left( x\right) =f\left(
x\right) -\sum_{j=1}^{n}\frac{f_{n-\alpha }^{\left( n-j\right) }\left(
a\right) }{\Gamma \left( \alpha -j+1\right) }\left( x-a\right) ^{\alpha -j}, 
\]%
holds almost everywhere on $\left[ a,b\right] $.

\bigskip

\textbf{Lemma 7. }[30]\textbf{\ (}Schaefer's fixed point theorem) Let $X$ be
a Banach space and $T:X\rightarrow X$ be a continuous and compact mapping.
If the set%
\[
\left\{ x\in X:x=\lambda T\left( x\right) \text{ for some }\lambda \in \left[
0,1\right] \right\} 
\]%
is bounded, then $T$ has a fixed point.

\bigskip

\textbf{3. Existence Result}

Recently, problems involving the fractional $p$-Laplacian operator have been
of great interest and this subject is studied by many mathematician. Taiyong
Chen, Wenbin Liu discuss the existence of solutions for the anti-periodic
boundary value problem of a fractional $p$-Laplacian equation given by%
\begin{eqnarray*}
D_{0,+}^{\beta }\phi _{p}\left( D_{0,+}^{\alpha }x\left( t\right) \right)
&=&f\left( t,x\left( t\right) \right) ,\text{ \ \ \ \ \ }t\in \left[ 0,1%
\right] , \\
x\left( 0\right) &=&-x\left( 1\right) ,\text{ \ \ \ \ \ }D_{0,+}^{\alpha
}x\left( 0\right) =-D_{0,+}^{\alpha }x\left( 1\right) ,
\end{eqnarray*}%
where $0<\alpha $, $\beta \leq 1$, $1<\alpha +\beta \leq 2$, $%
D_{0,+}^{\alpha }$ is a Caputo fractional derivative. We investigate
fractional $p$-Laplacian Sturm-Liouville problem having diffusion operator
with impulsive conditions and establish a theorem on existence of solution
the problem. We use Riemann-Liouville and Caputo fractional derivatives and
potential function.

\bigskip

\textbf{Theorem 8.} For a given $y\in PC\left( J,R\right) $. A function $y$
is a solution of fractional Sturm-Liouvlle problem%
\[
-D_{0,+}^{\beta }\phi _{p}\;^{C}D_{0,+}^{\alpha }y\left( t\right) +\left(
2\lambda p\left( t\right) +q\left( t\right) \right) y\left( t\right) =0, 
\]%
\[
\Delta y\mid _{t=t_{k}}=I_{k}\left( y\left( t_{k}\right) \right) ,\text{ }%
\Delta y^{\prime }\mid _{t=t_{k}}=I_{k}^{\ast }\left( y\left( t_{k}\right)
\right) ,\text{ }t_{k}\in \left( 0,\pi \right) ,\text{ }k=1,2,..,n, 
\]%
\[
y\left( 0\right) +y^{\prime }\left( 0\right) =0,\text{ \ }y\left( \pi
\right) +y^{\prime }\left( \pi \right) =0, 
\]%
if and only if $y$ is a solution of the fractional integral equation 
\begin{eqnarray*}
y\left( t\right) &=&\dint\limits_{t_{k}}^{t}\frac{\left( t-s\right) ^{\alpha
-1}}{\Gamma \left( \alpha \right) }\phi _{q}I_{0,+}^{\beta }\left( 2\lambda
p\left( s\right) +q\left( s\right) \right) y\left( s\right) ds \\
&&+\sum_{i=1}^{n}\dint\limits_{t_{i-1}}^{t_{i}}\left( \frac{\left(
t-t_{i}\right) \left( t_{i}-s\right) ^{\alpha -2}}{\Gamma \left( \alpha
-1\right) }+\frac{\left( t_{i}-s\right) ^{\alpha -1}}{\Gamma \left( \alpha
\right) }\right) \phi _{q}I_{0,+}^{\beta }\left( 2\lambda p\left( s\right)
+q\left( s\right) \right) y\left( s\right) ds \\
&&+\left[ \frac{\left( 1-t\right) }{\pi }\dint\limits_{t_{k}}^{\pi }\frac{%
\left( \pi -s\right) ^{\alpha -1}}{\Gamma \left( \alpha \right) }\phi
_{q}I_{0,+}^{\beta }\left( 2\lambda p\left( s\right) +q\left( s\right)
\right) y\left( s\right) ds\right]
\end{eqnarray*}

$%
\begin{array}{c}
+\left[ \frac{\left( 1-t\right) }{\pi }\sum\limits_{i=1}^{n}\dint%
\limits_{t_{i-1}}^{t_{i}}\left( \frac{\left( \pi -t_{i}\right) \left(
t_{i}-s\right) ^{\alpha -2}}{\Gamma \left( \alpha -1\right) }+\frac{\left(
t_{i}-s\right) ^{\alpha -1}}{\Gamma \left( \alpha \right) }\right) \phi
_{q}I_{0,+}^{\beta }\left( 2\lambda p\left( s\right) +q\left( s\right)
\right) y\left( s\right) ds\right]%
\end{array}%
$%
\begin{eqnarray}
&&+\left[ \frac{\left( 1-t\right) }{\pi }\dint\limits_{t_{k}}^{\pi }\frac{%
\left( \pi -s\right) ^{\alpha -2}}{\Gamma \left( \alpha -1\right) }\phi
_{q}I_{0,+}^{\beta }\left( 2\lambda p\left( s\right) +q\left( s\right)
\right) y\left( s\right) ds\right]  \nonumber \\
&&+\frac{\left( 1-t\right) }{\pi }\left[ \sum_{i=1}^{n}\dint%
\limits_{t_{i-1}}^{t_{i}}\frac{\left( t_{i}-s\right) ^{\alpha -2}}{\Gamma
\left( \alpha -1\right) }\phi _{q}I_{0,+}^{\beta }\left( 2\lambda p\left(
s\right) +q\left( s\right) \right) y\left( s\right) ds\right]  \nonumber \\
&&+\left[ \frac{\left( \pi +1-t\right) }{\pi }\sum_{i=1}^{n}I_{i}\left(
y\left( t_{i}\right) \right) +\frac{\left( \pi +1-t\right) }{\pi }%
\sum_{i=1}^{n}I_{i}^{\ast }\left( y\left( t_{i}\right) \right) \left(
1-t_{i}\right) \right] .  \TCItag{4}
\end{eqnarray}%
\textbf{Proof. }Assuming $y$ satisfies $\left( 1\right) -\left( 3\right) $.
Using Lemma 6 and Lemma 5, for some constants $b_{0},b_{1}\in 
\mathbb{R}
$, $t\in \left( 0,t_{1}\right] $ we have%
\begin{eqnarray}
y\left( t\right) &=&I_{0,+}^{\alpha }\phi _{q}I_{0,+}^{\beta }\left(
2\lambda p\left( t\right) +q\left( t\right) \right) y\left( t\right)
+b_{0}+b_{1}t,  \nonumber \\
&=&\dint\limits_{0}^{t}\frac{\left( t-s\right) ^{\alpha -1}}{\Gamma \left(
\alpha \right) }\phi _{q}I_{0,+}^{\beta }\left( 2\lambda p\left( s\right)
+q\left( s\right) \right) y\left( s\right) ds+b_{0}+b_{1}t,  \TCItag{5}
\end{eqnarray}%
It follows from $\left( 5\right) $ that 
\[
y^{\prime }\left( t\right) =\dint\limits_{0}^{t}\frac{\left( t-s\right)
^{\alpha -2}}{\Gamma \left( \alpha -1\right) }\phi _{q}I_{0,+}^{\beta
}\left( 2\lambda p\left( s\right) +q\left( s\right) \right) y\left( s\right)
ds+b_{1}, 
\]%
if $t\in \left( t_{1},t_{2}\right] $ and $c_{0}$,$c_{1}\in 
\mathbb{R}
$ are arbitrary constants then we have%
\begin{eqnarray*}
y\left( t\right) &=&\dint\limits_{t_{1}}^{t}\frac{\left( t-s\right) ^{\alpha
-1}}{\Gamma \left( \alpha \right) }\phi _{q}I_{0,+}^{\beta }\left( 2\lambda
p\left( s\right) +q\left( s\right) \right) y\left( s\right)
ds+c_{0}+c_{1}\left( t-t_{1}\right) , \\
y^{\prime }\left( t\right) &=&\dint\limits_{t_{1}}^{t}\frac{\left(
t-s\right) ^{\alpha -2}}{\Gamma \left( \alpha -1\right) }\phi
_{q}I_{0,+}^{\beta }\left( 2\lambda p\left( s\right) +q\left( s\right)
\right) y\left( s\right) ds+c_{1},
\end{eqnarray*}%
using the impulse conditions $(2)$%
\begin{eqnarray*}
c_{0} &=&\dint\limits_{0}^{t_{1}}\frac{\left( t_{1}-s\right) ^{\alpha -1}}{%
\Gamma \left( \alpha \right) }\phi _{q}I_{0,+}^{\beta }\left( 2\lambda
p\left( s\right) +q\left( s\right) \right) y\left( s\right) ds+b_{0} \\
&&+b_{1}t_{1}+I_{1}\left( y\left( t_{1}\right) \right) ,
\end{eqnarray*}%
\[
c_{1}=\dint\limits_{0}^{t_{1}}\frac{\left( t_{1}-s\right) ^{\alpha -2}}{%
\Gamma \left( \alpha -1\right) }\phi _{q}I_{0,+}^{\beta }\left( 2\lambda
p\left( s\right) +q\left( s\right) \right) y\left( s\right)
ds+b_{1}+I_{1}^{\ast }\left( y\left( t_{1}\right) \right) , 
\]%
thus,%
\begin{eqnarray*}
y\left( t\right) &=&\dint\limits_{t_{1}}^{t}\frac{\left( t-s\right) ^{\alpha
-1}}{\Gamma \left( \alpha \right) }\phi _{q}I_{0,+}^{\beta }\left( 2\lambda
p\left( s\right) +q\left( s\right) \right) y\left( s\right) ds \\
&&+\dint\limits_{0}^{t_{1}}\left( \frac{\left( t-t_{1}\right) \left(
t_{1}-s\right) ^{\alpha -2}}{\Gamma \left( \alpha -1\right) }+\frac{\left(
t_{1}-s\right) ^{\alpha -1}}{\Gamma \left( \alpha \right) }\right) \phi
_{q}I_{0,+}^{\beta }\left( 2\lambda p\left( s\right) +q\left( s\right)
\right) y\left( s\right) ds \\
&&+b_{0}+b_{1}t+I_{1}\left( y\left( t_{1}\right) \right) +I_{1}^{\ast
}\left( y\left( t_{1}\right) \right) \left( t-t_{1}\right) ,
\end{eqnarray*}%
repeating the process in this way, for $t\in \left( t_{k},t_{k+1}\right] $,
we have 
\begin{eqnarray}
y\left( t\right) &=&\dint\limits_{t_{k}}^{t}\left[ \frac{\left( t-s\right)
^{\alpha -1}}{\Gamma \left( \alpha \right) }\phi _{q}I_{0,+}^{\beta }\left(
2\lambda p\left( s\right) +q\left( s\right) \right) y\left( s\right) \right]
ds  \nonumber \\
&&+\left[ \sum_{i=1}^{n}\dint\limits_{t_{i-1}}^{t_{i}}\left( \frac{\left(
t-t_{i}\right) \left( t_{i}-s\right) ^{\alpha -2}}{\Gamma \left( \alpha
-1\right) }+\frac{\left( t_{i}-s\right) ^{\alpha -1}}{\Gamma \left( \alpha
\right) }\right) \phi _{q}I_{0,+}^{\beta }\left( 2\lambda p\left( s\right)
+q\left( s\right) \right) y\left( s\right) ds\right]  \nonumber \\
&&+b_{0}+b_{1}t+\sum_{i=1}^{n}I_{i}\left( y\left( t_{i}\right) \right)
+\sum_{i=1}^{n}I_{i}^{\ast }\left( y\left( t_{i}\right) \right) \left(
t-t_{i}\right) ,  \TCItag{6}
\end{eqnarray}%
and%
\begin{eqnarray*}
y^{\prime }\left( t\right) &=&\left[ \dint\limits_{t_{k}}^{t}\frac{\left(
t-s\right) ^{\alpha -2}}{\Gamma \left( \alpha -1\right) }\phi
_{q}I_{0,+}^{\beta }\left( 2\lambda p\left( s\right) +q\left( s\right)
\right) y\left( s\right) ds\right] \\
&&+\left[ \sum_{i=1}^{n}\dint\limits_{t_{i-1}}^{t_{i}}\frac{\left(
t_{i}-s\right) ^{\alpha -2}}{\Gamma \left( \alpha -1\right) }\phi
_{q}I_{0,+}^{\beta }\left( 2\lambda p\left( s\right) +q\left( s\right)
\right) y\left( s\right) ds\right] \\
&&+b_{1}+\sum_{i=1}^{n}I_{i}^{\ast }\left( y\left( t_{i}\right) \right) ,
\end{eqnarray*}%
applying the boundary condition $y\left( 0\right) +y^{\prime }\left(
0\right) =0,$ \ $y\left( \pi \right) +y^{\prime }\left( \pi \right) =0$, we
find that \newpage 
\begin{eqnarray}
&&b_{1}=-\frac{1}{\pi }\dint\limits_{t_{k}}^{\pi }\frac{\left( \pi -s\right)
^{\alpha -1}}{\Gamma \left( \alpha \right) }\phi _{q}I_{0,+}^{\beta }\left(
2\lambda p\left( s\right) +q\left( s\right) \right) y\left( s\right) ds 
\nonumber \\
&&-\left[ \frac{1}{\pi }\sum_{i=1}^{n}\dint\limits_{t_{i-1}}^{t_{i}}\left( 
\frac{\left( \pi -t_{i}\right) \left( t_{i}-s\right) ^{\alpha -2}}{\Gamma
\left( \alpha -1\right) }+\frac{\left( t_{i}-s\right) ^{\alpha -1}}{\Gamma
\left( \alpha \right) }\right) \phi _{q}I_{0,+}^{\beta }\left( 2\lambda
p\left( s\right) +q\left( s\right) \right) y\left( s\right) ds\right] \\
&&-\left[ \frac{1}{\pi }\sum_{i=1}^{n}I_{i}\left( y\left( t_{i}\right)
\right) -\frac{1}{\pi }\sum_{i=1}^{n}I_{i}^{\ast }\left( y\left(
t_{i}\right) \right) \left( \pi -t_{i}\right) \right]  \nonumber \\
&&-\left[ \frac{1}{\pi }\dint\limits_{t_{k}}^{\pi }\frac{\left( \pi
-s\right) ^{\alpha -2}}{\Gamma \left( \alpha -1\right) }\phi
_{q}I_{0,+}^{\beta }\left( 2\lambda p\left( s\right) +q\left( s\right)
\right) y\left( s\right) ds\right]  \nonumber \\
&&-\left[ \frac{1}{\pi }\sum_{i=1}^{n}\dint\limits_{t_{i-1}}^{t_{i}}\frac{%
\left( t_{i}-s\right) ^{\alpha -2}}{\Gamma \left( \alpha -1\right) }\phi
_{q}I_{0,+}^{\beta }\left( 2\lambda p\left( s\right) +q\left( s\right)
\right) y\left( s\right) ds-\frac{1}{\pi }\sum_{i=1}^{n}I_{i}^{\ast }\left(
y\left( t_{i}\right) \right) \right] ,  \TCItag{7}
\end{eqnarray}%
\[
b_{0}=\frac{1}{\pi }\dint\limits_{t_{k}}^{\pi }\frac{\left( \pi -s\right)
^{\alpha -1}}{\Gamma \left( \alpha \right) }\phi _{q}I_{0,+}^{\beta }\left(
2\lambda p\left( s\right) +q\left( s\right) \right) y\left( s\right) ds 
\]%
\begin{eqnarray*}
&&+\left[ \frac{1}{\pi }\sum_{i=1}^{n}\dint\limits_{t_{i-1}}^{t_{i}}\left( 
\frac{\left( \pi -t_{i}\right) \left( t_{i}-s\right) ^{\alpha -2}}{\Gamma
\left( \alpha -1\right) }+\frac{\left( t_{i}-s\right) ^{\alpha -1}}{\Gamma
\left( \alpha \right) }\right) \phi _{q}I_{0,+}^{\beta }\left( 2\lambda
p\left( s\right) +q\left( s\right) \right) y\left( s\right) ds\right] \\
&&+\left[ \frac{1}{\pi }\sum_{i=1}^{n}I_{i}\left( y\left( t_{i}\right)
\right) +\frac{1}{\pi }\sum_{i=1}^{n}I_{i}^{\ast }\left( y\left(
t_{i}\right) \right) \left( \pi -t_{i}\right) \right]
\end{eqnarray*}%
\begin{eqnarray}
&&+\frac{1}{\pi }\dint\limits_{t_{k}}^{\pi }\frac{\left( \pi -s\right)
^{\alpha -2}}{\Gamma \left( \alpha -1\right) }\phi _{q}I_{0,+}^{\beta
}\left( 2\lambda p\left( s\right) +q\left( s\right) \right) y\left( s\right)
ds  \nonumber \\
&&+\frac{1}{\pi }\left[ \sum_{i=1}^{n}\dint\limits_{t_{i-1}}^{t_{i}}\frac{%
\left( t_{i}-s\right) ^{\alpha -2}}{\Gamma \left( \alpha -1\right) }\phi
_{q}I_{0,+}^{\beta }\left( 2\lambda p\left( s\right) +q\left( s\right)
\right) y\left( s\right) ds+\frac{1}{\pi }\sum_{i=1}^{n}I_{i}^{\ast }\left(
y\left( t_{i}\right) \right) \right] ,  \TCItag{8}
\end{eqnarray}%
substituting $(7),(8)$ into $(6)$, we obtain $(4)$. The proof completes.

\bigskip

\textbf{Theorem 9. }Presume that

\begin{description}
\item[(H$_{\text{1}}$)] There exist constants $N,R,M>0$ such that%
\[
\text{ }\left\vert \lambda \right\vert \leq N,\text{ }\left\vert p\left(
t\right) \right\vert \leq R,\text{ }\left\vert q\left( t\right) \right\vert
\leq M\text{ for each }t\in J. 
\]

\item[(H$_{\text{2}}$)] The functions $I_{k},I_{k}^{\ast }:%
\mathbb{R}
\longrightarrow 
\mathbb{R}
$ are continuous and there exist constant $r_{1},r_{2}>0$ such that$%
\left\vert I_{k}\left( y\right) \right\vert <r_{1},$ $\left\vert I_{k}^{\ast
}\left( y\right) \right\vert <r_{2},$ $k=1,...,n,\left\vert I_{k}\left(
y\right) \right\vert <r_{1},\left\vert I_{k}^{\ast }\left( y\right)
\right\vert <r_{2},k=1,...,n,$%
\[
\left\vert I_{k}\left( y\right) \right\vert <r_{1},\left\vert I_{k}^{\ast
}\left( y\right) \right\vert <r_{2},k=1,...,n, 
\]
\end{description}

\noindent then the $\left( 1\right) -\left( 3\right) $ problem has at least
one solution on $J.$

\textbf{Proof. }Define the operator $T:PC\left( J,R\right) \rightarrow
PC\left( J,R\right) $ as%
\begin{eqnarray*}
&&T\left( y\left( t\right) \right) 
\begin{array}{c}
=%
\end{array}%
\dint\limits_{t_{k}}^{t}\frac{\left( t-s\right) ^{\alpha -1}}{\Gamma \left(
\alpha \right) }\phi _{q}I_{0,+}^{\beta }\left( 2\lambda p\left( s\right)
+q\left( s\right) \right) y\left( s\right) ds \\
&&+\left[ \sum_{i=1}^{n}\dint\limits_{t_{i-1}}^{t_{i}}\left( \frac{\left(
t-t_{i}\right) \left( t_{i}-s\right) ^{\alpha -2}}{\Gamma \left( \alpha
-1\right) }+\frac{\left( t_{i}-s\right) ^{\alpha -1}}{\Gamma \left( \alpha
\right) }\right) \phi _{q}I_{0,+}^{\beta }\left( 2\lambda p\left( s\right)
+q\left( s\right) \right) y\left( s\right) ds\right] \\
&&+\left[ \frac{\left( 1-t\right) }{\pi }\dint\limits_{t_{k}}^{\pi }\frac{%
\left( \pi -s\right) ^{\alpha -1}}{\Gamma \left( \alpha \right) }\phi
_{q}I_{0,+}^{\beta }\left( 2\lambda p\left( s\right) +q\left( s\right)
\right) y\left( s\right) ds\right] \\
&&+\left[ \frac{\left( 1-t\right) }{\pi }\tsum_{i=1}^{n}\tint%
\limits_{t_{i-1}}^{t_{i}}\left( \frac{\left( \pi -t_{i}\right) \left(
t_{i}-s\right) ^{\alpha -2}}{\Gamma \left( \alpha -1\right) }+\frac{\left(
t_{i}-s\right) ^{\alpha -1}}{\Gamma \left( \alpha \right) }\right) \phi
_{q}I_{0,+}^{\beta }\left( 2\lambda p\left( s\right) +q\left( s\right)
\right) y\left( s\right) ds\right] \\
&&+\left[ \frac{\left( 1-t\right) }{\pi }\dint\limits_{t_{k}}^{\pi }\frac{%
\left( \pi -s\right) ^{\alpha -2}}{\Gamma \left( \alpha -1\right) }\phi
_{q}I_{0,+}^{\beta }\left( 2\lambda p\left( s\right) +q\left( s\right)
\right) y\left( s\right) ds\right] \\
&&+\left[ \frac{\left( 1-t\right) }{\pi }\sum_{i=1}^{n}\dint%
\limits_{t_{i-1}}^{t_{i}}\frac{\left( t_{i}-s\right) ^{\alpha -2}}{\Gamma
\left( \alpha -1\right) }\phi _{q}I_{0,+}^{\beta }\left( 2\lambda p\left(
s\right) +q\left( s\right) \right) y\left( s\right) ds\right] \\
&&+\left[ \frac{\left( \pi +1-t\right) }{\pi }\sum_{i=1}^{n}I_{i}\left(
y\left( t_{i}\right) \right) +\frac{\left( \pi +1-t\right) }{\pi }%
\sum_{i=1}^{n}I_{i}^{\ast }\left( y\left( t_{i}\right) \right) \left(
1-t_{i}\right) \right] .
\end{eqnarray*}%
Now, to prove that $T$ has a fixed point, we use Schaefer fixed point
theorem and it will be proven in four steps.

\begin{description}
\item[Step 1:] $T$\textbf{\ }is continuous.
\end{description}

\noindent Let $\left\{ y_{n}\right\} $ be a sequence such that $%
y_{n}\rightarrow y$ in $PC\left( J,R\right) $. Then for each $t\in J$%
\[
\left\vert T\left( y_{n}\right) \left( t\right) -T\left( y\right) \left(
t\right) \right\vert \leq \dint\limits_{t_{k}}^{t}\frac{\left( t-s\right)
^{\alpha -1}}{\Gamma \left( \alpha \right) }\phi _{q}\left\vert
I_{0,+}^{\beta }\left( 2\lambda p\left( s\right) +q\left( s\right) \right)
\left( y_{n}\left( s\right) -y\left( s\right) \right) \right\vert ds 
\]%
$%
\begin{array}{c}
+\sum\limits_{i=1}^{n}\dint\limits_{t_{i-1}}^{t_{i}}\left( \frac{\left(
t-t_{i}\right) \left( t_{i}-s\right) ^{\alpha -2}}{\Gamma \left( \alpha
-1\right) }+\frac{\left( t_{i}-s\right) ^{\alpha -1}}{\Gamma \left( \alpha
\right) }\right) \phi _{q}\left\vert I_{0,+}^{\beta }\left( 2\lambda p\left(
s\right) +q\left( s\right) \right) \left( y_{n}\left( s\right) -y\left(
s\right) \right) \right\vert ds \\ 
+\frac{\left( 1-t\right) }{\pi }\dint\limits_{t_{k}}^{\pi }\left( \frac{%
\left( \pi -s\right) ^{\alpha -1}}{\Gamma \left( \alpha \right) }+\frac{%
\left( \pi -s\right) ^{\alpha -2}}{\Gamma \left( \alpha -1\right) }\right)
\phi _{q}\left\vert I_{0,+}^{\beta }\left( 2\lambda p\left( s\right)
+q\left( s\right) \right) \left( y_{n}\left( s\right) -y\left( s\right)
\right) \right\vert ds \\ 
+\frac{\left( 1-t\right) }{\pi }\sum\limits_{i=1}^{n}\dint%
\limits_{t_{i-1}}^{t_{i}}\left( \frac{\left( \pi -t_{i}\right) \left(
t_{i}-s\right) ^{\alpha -2}}{\Gamma \left( \alpha -1\right) }+\frac{\left(
t_{i}-s\right) ^{\alpha -1}}{\Gamma \left( \alpha \right) }\right) \phi
_{q}\left\vert I_{0,+}^{\beta }\left( 2\lambda p\left( s\right) +q\left(
s\right) \right) \left( y_{n}\left( s\right) -y\left( s\right) \right)
\right\vert ds \\ 
+\frac{\left( 1-t\right) }{\pi }\sum\limits_{i=1}^{n}\dint%
\limits_{t_{i-1}}^{t_{i}}\frac{\left( t_{i}-s\right) ^{\alpha -2}}{\Gamma
\left( \alpha -1\right) }\phi _{q}\left\vert I_{0,+}^{\beta }\left( 2\lambda
p\left( s\right) +q\left( s\right) \right) \left( y_{n}\left( s\right)
-y\left( s\right) \right) \right\vert ds \\ 
+\frac{\left( \pi +1-t\right) }{\pi }\sum\limits_{i=1}^{n}\left\vert
I_{i}\left( y_{n}\left( t_{i}\right) \right) -I_{i}\left( y\left(
t_{i}\right) \right) \right\vert +\frac{\left( \pi +1-t\right) }{\pi }%
\sum\limits_{i=1}^{n}\left\vert I_{i}^{\ast }\left( y_{n}\left( t_{i}\right)
\right) -I_{i}^{\ast }\left( y\left( t_{i}\right) \right) \right\vert \left(
1-t_{i}\right) ,%
\end{array}%
$

Since $I_{k},I_{k}^{\ast }$, $k=1,...,n,$ are continuous functions, we have%
\[
\left\Vert T\left( y_{n}\right) -T\left( y\right) \right\Vert _{\infty
}\rightarrow 0,\text{ as }n\rightarrow \infty . 
\]

\begin{description}
\item[Step 2:] $T$ operator bounded on bounded sets of $PC\left( J,R\right)
. $
\end{description}

In fact, it is enough to show that for any $\nu >0$ there exists a positive
constant $\delta $ such that for each $y\in B=\left\{ y\in
PC(J,R):\left\Vert y\right\Vert _{\infty }<\nu \right\} $ we have $%
\left\Vert T\left( y\right) \right\Vert _{\infty }\leq \delta $. There
exists constant $K>0$ such that $\left\vert I_{0,+}^{\beta }\left( 2\lambda
p\left( s\right) +q\left( s\right) \right) y\left( s\right) \right\vert \leq
K$. By (H$_{\text{1}}$) and (H$_{\text{2}}$), we have for each $t\in J$

$%
\begin{array}{c}
\left\vert Ty\left( t\right) \right\vert \leq \left\vert
\dint\limits_{t_{k}}^{t}\frac{\left( t-s\right) ^{\alpha -1}}{\Gamma \left(
\alpha \right) }\phi _{q}I_{0,+}^{\beta }\left( 2\lambda p\left( s\right)
+q\left( s\right) \right) y\left( s\right) ds\right\vert%
\end{array}%
$%
\begin{eqnarray*}
&&+\left[ \left\vert
\sum\limits_{i=1}^{n}\dint\limits_{t_{i-1}}^{t_{i}}\left( \frac{\left(
t-t_{i}\right) \left( t_{i}-s\right) ^{\alpha -2}}{\Gamma \left( \alpha
-1\right) }+\frac{\left( t_{i}-s\right) ^{\alpha -1}}{\Gamma \left( \alpha
\right) }\right) \phi _{q}I_{0,+}^{\beta }\left( 2\lambda p\left( s\right)
+q\left( s\right) \right) y\left( s\right) ds\right\vert \right] \\
&&+\left[ \left\vert \frac{\left( 1-t\right) }{\pi }\dint\limits_{t_{n}}^{%
\pi }\left( \frac{\left( \pi -s\right) ^{\alpha -1}}{\Gamma \left( \alpha
\right) }+\frac{\left( \pi -s\right) ^{\alpha -2}}{\Gamma \left( \alpha
-1\right) }\right) \phi _{q}I_{0,+}^{\beta }\left( 2\lambda p\left( s\right)
+q\left( s\right) \right) y\left( s\right) ds\right\vert \right] \\
&&+\left[ \left\vert \frac{\left( 1-t\right) }{\pi }\sum\limits_{i=1}^{n}%
\dint\limits_{t_{i-1}}^{t_{i}}\left( \frac{\left( \pi -t_{i}\right) \left(
t_{i}-s\right) ^{\alpha -2}}{\Gamma \left( \alpha -1\right) }+\frac{\left(
t_{i}-s\right) ^{\alpha -1}}{\Gamma \left( \alpha \right) }\right) \phi
_{q}I_{0,+}^{\beta }\left( 2\lambda p\left( s\right) +q\left( s\right)
\right) y\left( s\right) ds\right\vert \right] \\
&&+\left[ \left\vert \frac{\left( 1-t\right) }{\pi }\sum\limits_{i=1}^{n}%
\dint\limits_{t_{i-1}}^{t_{i}}\frac{\left( t_{i}-s\right) ^{\alpha -2}}{%
\Gamma \left( \alpha -1\right) }\phi _{q}I_{0,+}^{\beta }\left( 2\lambda
p\left( s\right) +q\left( s\right) \right) y\left( s\right) ds\right\vert %
\right] \\
&&+\left[ \left\vert \frac{\pi +1-t}{\pi }\sum_{i=1}^{n}I_{i}\left( y\left(
t_{i}\right) \right) \right\vert +\left\vert \frac{\pi +1-t}{\pi }%
\sum_{i=1}^{n}I_{i}^{\ast }\left( y\left( t_{i}\right) \right) \left(
1-t_{i}\right) \right\vert \right] \\
&\leq &\left[ \frac{K^{q-1}}{\Gamma \left( \alpha \right) }%
\dint\limits_{t_{n}}^{t}\left( t-s\right) ^{\alpha
-1}ds+K^{q-1}\sum_{i=1}^{n}\dint\limits_{t_{i-1}}^{t_{i}}\left( \frac{\left(
t-t_{i}\right) \left( t_{i}-s\right) ^{\alpha -2}}{\Gamma \left( \alpha
-1\right) }+\frac{\left( t_{i}-s\right) ^{\alpha -1}}{\Gamma \left( \alpha
\right) }\right) ds\right] \\
&&+\left[ \frac{1}{\pi }K^{q-1}\dint\limits_{t_{n}}^{\pi }\left( \frac{%
\left( \pi -s\right) ^{\alpha -1}}{\Gamma \left( \alpha \right) }+\frac{%
\left( \pi -s\right) ^{\alpha -2}}{\Gamma \left( \alpha -1\right) }\right) ds%
\right]
\end{eqnarray*}%
\begin{eqnarray*}
&&+\frac{1}{\pi }K^{q-1}\sum_{i=1}^{n}\dint\limits_{t_{i-1}}^{t_{i}}\left( 
\frac{\left( \pi -t_{i}\right) \left( t_{i}-s\right) ^{\alpha -2}}{\Gamma
\left( \alpha -1\right) }+\frac{\left( t_{i}-s\right) ^{\alpha -1}}{\Gamma
\left( \alpha \right) }\right) ds \\
&&+\frac{1}{\pi }K^{q-1}\sum_{i=1}^{n}\dint\limits_{t_{i-1}}^{t_{i}}\frac{%
\left( t_{i}-s\right) ^{\alpha -2}}{\Gamma \left( \alpha -1\right) }ds+\frac{%
\left( 1+\pi \right) nr_{1}}{\pi }+\frac{\left( 1+\pi \right) nr_{2}}{\pi }
\\
&\leq &\frac{K^{q-1}\pi ^{\alpha }}{\Gamma \left( \alpha +1\right) }+\frac{%
nK^{q-1}\pi \pi ^{\alpha -1}}{\Gamma \left( \alpha \right) }+\frac{%
nK^{q-1}\pi ^{\alpha }}{\Gamma \left( \alpha +1\right) }+\frac{K^{q-1}\pi
^{\alpha }}{\Gamma \left( \alpha +1\right) \pi }
\end{eqnarray*}

\begin{eqnarray*}
&&+\frac{K^{q-1}\pi ^{\alpha -1}}{\Gamma \left( \alpha \right) \pi }+\frac{%
nK^{q-1}\pi \pi ^{\alpha -1}}{\Gamma \left( \alpha \right) \pi }+\frac{%
nK^{q-1}\pi ^{\alpha }}{\Gamma \left( \alpha +1\right) \pi }+\frac{%
nK^{q-1}\pi ^{\alpha -1}}{\Gamma \left( \alpha \right) \pi } \\
&&+\frac{\left( 1+\pi \right) nr_{1}}{\pi }+\frac{\left( 1+\pi \right) nr_{2}%
}{\pi } \\
&\leq &K^{q-1}\pi ^{\alpha }\left[ \frac{\left( n+1\right) \left( \pi
+1\right) }{\pi \Gamma \left( \alpha +1\right) }+\frac{n\left( \pi +1\right) 
}{\pi \Gamma \left( \alpha \right) }\right] +K^{q-1}\pi ^{\alpha -1}\frac{%
\left( n+1\right) }{\Gamma \left( \alpha \right) \pi } \\
&&+\frac{n\left( \pi +1\right) \left( r_{1}+r_{2}\right) }{\pi }
\end{eqnarray*}%
Thus%
\begin{eqnarray*}
\left\Vert T\left( y\right) \right\Vert _{\infty } &\leq &K^{q-1}\pi
^{\alpha }\left[ \frac{\left( n+1\right) \left( \pi +1\right) }{\pi \Gamma
\left( \alpha +1\right) }+\frac{n\left( \pi +1\right) }{\pi \Gamma \left(
\alpha \right) }\right] +K^{q-1}\pi ^{\alpha -1}\frac{\left( n+1\right) }{%
\Gamma \left( \alpha \right) \pi } \\
&&+\frac{n\left( \pi +1\right) \left( r_{1}+r_{2}\right) }{\pi }%
\begin{array}{c}
=%
\end{array}%
\delta .
\end{eqnarray*}

\begin{description}
\item[Step 3] $T$ operator bounded equicontinuous on sets of $PC\left(
J,R\right) .$
\end{description}

Let $\tau _{1},\tau _{2}\in J,$ $\tau _{1}<\tau _{2}$, $B$ be a bounded set
of $PC\left( J,R\right) $ as in Step 2, and let $y\in B$. Then

$%
\begin{array}{c}
\left\vert T\left( y\right) \left( \tau _{2}\right) -T\left( y\right) \left(
\tau _{1}\right) \right\vert \leq \left[ \dint\limits_{t_{k}}^{\tau _{1}}%
\frac{\left( \tau _{2}-s\right) ^{\alpha -1}}{\Gamma \left( \alpha \right) }-%
\frac{\left( \tau _{1}-s\right) ^{\alpha -1}}{\Gamma \left( \alpha \right) }%
\phi _{q}\left\vert I_{0,+}^{\beta }\left( 2\lambda p\left( s\right)
+q\left( s\right) \right) y\left( s\right) \right\vert ds\right] \\ 
+\left[ \dint\limits_{\tau _{1}}^{\tau _{2}}\frac{\left( \tau _{2}-s\right)
^{\alpha -1}}{\Gamma \left( \alpha \right) }\phi _{q}\left\vert
I_{0,+}^{\beta }\left( 2\lambda p\left( s\right) +q\left( s\right) \right)
y\left( s\right) \right\vert ds\right] \\ 
+\left[ \frac{\left( \tau _{2}-\tau _{1}\right) }{\pi }\sum\limits_{i=1}^{n}%
\dint\limits_{t_{i-1}}^{t_{i}}\frac{\left( t_{i}-s\right) ^{\alpha -2}}{%
\Gamma \left( \alpha -1\right) }\phi _{q}\left\vert I_{0,+}^{\beta }\left(
2\lambda p\left( s\right) +q\left( s\right) \right) y\left( s\right)
\right\vert ds\right] \\ 
+\left[ \frac{\left( \tau _{2}-\tau _{1}\right) }{\pi }\dint%
\limits_{t_{k}}^{\pi }\frac{\left( \pi -s\right) ^{\alpha -1}}{\Gamma \left(
\alpha \right) }+\frac{\left( \pi -s\right) ^{\alpha -2}}{\Gamma \left(
\alpha -1\right) }\phi _{q}\left\vert I_{0,+}^{\beta }\left( 2\lambda
p\left( s\right) +q\left( s\right) \right) y\left( s\right) \right\vert ds%
\right] \\ 
+\left[ \frac{\left( \tau _{2}-\tau _{1}\right) }{\pi }\sum\limits_{i=1}^{n}%
\dint\limits_{t_{i-1}}^{t_{i}}\left( \frac{\left( \pi -t_{i}\right) \left(
t_{i}-s\right) ^{\alpha -2}}{\Gamma \left( \alpha -1\right) }+\frac{\left(
t_{i}-s\right) ^{\alpha -1}}{\Gamma \left( \alpha \right) }\right) \phi
_{q}\left\vert I_{0,+}^{\beta }\left( 2\lambda p\left( s\right) +q\left(
s\right) \right) y\left( s\right) \right\vert ds\right] \\ 
+\left[ \frac{\left( \tau _{2}-\tau _{1}\right) }{\pi }\sum\limits_{i=1}^{n}%
\dint\limits_{t_{i-1}}^{t_{i}}\frac{\left( t_{i}-s\right) ^{\alpha -2}}{%
\Gamma \left( \alpha -1\right) }\phi _{q}\left\vert I_{0,+}^{\beta }\left(
2\lambda p\left( s\right) +q\left( s\right) \right) y\left( s\right)
\right\vert ds\right] \\ 
+\left[ \frac{\left( \tau _{2}-\tau _{1}\right) }{\pi }\sum\limits_{i=1}^{n}%
\left\vert I_{i}\left( y\left( t_{i}\right) \right) \right\vert +\frac{%
\left( \tau _{2}-\tau _{1}\right) }{\pi }\sum\limits_{i=1}^{n}\left\vert
I_{i}^{\ast }\left( y\left( t_{i}\right) \right) \right\vert \left(
1-t_{i}\right) \right] ,%
\end{array}%
$

As a result of Step 1 to Step 3 by Arzela-Ascoli theorem, we can deduce that 
$T:PC\left( J,R\right) \rightarrow PC\left( J,R\right) $ is continuous and
completely continuous. Furthermore, as $\tau _{2}\rightarrow \tau _{1}$, $T$
operator is equicontinuous.

\begin{description}
\item[Step 4:] Now, let's show that the set
\end{description}

\[
L=\left \{ y\in PC\left[ J,R\right] :y=\theta T\left( y\right) ,\text{ }%
0<\theta <1\right \} , 
\]%
is bounded.

Let $y\in L$. Then $y=\theta T\left( y\right) $, for some $0<\theta <1$.
Thus for each $t\in J,$ we have

$%
\begin{array}{c}
y\left( t\right) =\left[ \theta \dint\limits_{t_{n}}^{t}\frac{\left(
t-s\right) ^{\alpha -1}}{\Gamma \left( \alpha \right) }\phi
_{q}I_{0,+}^{\beta }\left( 2\lambda p\left( s\right) +q\left( s\right)
\right) y\left( s\right) ds\right] \\ 
+\left[ \theta \sum\limits_{i=1}^{n}\dint\limits_{t_{i-1}}^{t_{i}}\left( 
\frac{\left( t-t_{i}\right) \left( t_{i}-s\right) ^{\alpha -2}}{\Gamma
\left( \alpha -1\right) }+\frac{\left( t_{i}-s\right) ^{\alpha -1}}{\Gamma
\left( \alpha \right) }\right) \phi _{q}I_{0,+}^{\beta }\left( 2\lambda
p\left( s\right) +q\left( s\right) \right) y\left( s\right) ds\right] \\ 
+\left[ \frac{\left( 1-t\right) \theta }{\pi }\dint\limits_{t_{n}}^{\pi }%
\frac{\left( \pi -s\right) ^{\alpha -1}}{\Gamma \left( \alpha \right) }\phi
_{q}I_{0,+}^{\beta }\left( 2\lambda p\left( s\right) +q\left( s\right)
\right) y\left( s\right) ds\right] \\ 
+\left[ \frac{\left( 1-t\right) \theta }{\pi }\sum\limits_{i=1}^{n}\tint%
\limits_{t_{i-1}}^{t_{i}}\left( \frac{\left( \pi -t_{i}\right) \left(
t_{i}-s\right) ^{\alpha -2}}{\Gamma \left( \alpha -1\right) }+\frac{\left(
t_{i}-s\right) ^{\alpha -1}}{\Gamma \left( \alpha \right) }\right) \phi
_{q}I_{0,+}^{\beta }\left( 2\lambda p\left( s\right) +q\left( s\right)
\right) y\left( s\right) ds\right] \\ 
+\left[ \frac{\left( 1-t\right) \theta }{\pi }\dint\limits_{t_{n}}^{\pi }%
\frac{\left( \pi -s\right) ^{\alpha -2}}{\Gamma \left( \alpha -1\right) }%
\phi _{q}I_{0,+}^{\beta }\left( 2\lambda p\left( s\right) +q\left( s\right)
\right) y\left( s\right) ds\right]%
\end{array}%
$

$%
\begin{array}{c}
+\left[ \frac{\left( 1-t\right) \theta }{\pi }\sum\limits_{i=1}^{n}\dint%
\limits_{t_{i-1}}^{t_{i}}\frac{\left( t_{i}-s\right) ^{\alpha -2}}{\Gamma
\left( \alpha -1\right) }\phi _{q}I_{0,+}^{\beta }\left( 2\lambda p\left(
s\right) +q\left( s\right) \right) y\left( s\right) d\right] s%
\end{array}%
$%
\[
+\left[ \frac{\left( \pi +1-t\right) \theta }{\pi }\sum_{i=1}^{n}I_{i}\left(
y\left( t_{i}\right) \right) +\frac{\left( \pi +1-t\right) \theta }{\pi }%
\sum_{i=1}^{n}I_{i}^{\ast }\left( y\left( t_{i}\right) \right) \left(
1-t_{i}\right) \right] 
\]%
This implies by (H$_{\text{1}}$) and (H$_{\text{2}}$) (as in Step 2) that
for each $t\in J$ we have%
\begin{eqnarray*}
\left\vert y\left( t\right) \right\vert &\leq &K^{q-1}\pi ^{\alpha }\left[ 
\frac{\left( n+1\right) \left( \pi +1\right) }{\pi \Gamma \left( \alpha
+1\right) }+\frac{n\left( \pi +1\right) }{\pi \Gamma \left( \alpha \right) }%
\right] +K^{q-1}\pi ^{\alpha -1}\frac{\left( n+1\right) }{\Gamma \left(
\alpha \right) \pi } \\
&&+\frac{n\left( \pi +1\right) \left( r_{1}+r_{2}\right) }{\pi } \\
&=&\delta ,
\end{eqnarray*}%
Furthermore, the set $L$ is bounded. We conclude that $T$ has a fixed point
in the solution of the problem $\left( 1\right) -\left( 3\right) $,
according to Schafer's

\textbf{Conclusions}

In this paper, we investigate fractional $p$-Laplacian Sturm-Liouville
problem having diffusion operator with impulsive conditions at $\alpha \in
\left( 1,2\right] $. The derivatives are described in the Riemann-Liouville
and Caputo sense.The fractional impulsive differential equation and boundary
value problem involving fractional $p$-Laplacian is analyzed for the case of
our fractional Sturm-Liouville problem. This paper is dealt with
Sturm-Liouville problem involving impulsive differential equation of
fractional order. We show an explicit representation of solution of the
problem. By using Schaefer fixed point theorem we proved existence of
solution for fractional $p$-Laplacian Sturm-Liouville problem having
diffusion operator with impulsive conditions. We hope that our study will
make a new research in the area of fractional Stum-Liouville problems begin
with different boundary condition and many of its variations.

\bigskip

\textbf{References}

\begin{enumerate}
\item Baleanu, D., Machado, JAT, Luo, A., Fractional Dynamics and Control,
Springer, Berlin, 2012.

\item Sabatier, J., Agrawal, O.P., Machado, JAT, Advances in Fractional
Calculus: Theoretical Developments and Applications in Physics and
Engineering, Springer, Dordrecht, 2007.

\item Lakshmikantham, V., Leela, S., Vasundhara Devi, J., Theory of
Fractional Dynamic Systems, Cambridge Scientific Publishers, Cambridge, 2009.

\item Kilbas, A.A., Srivastava, H.M. and Trujillo, J.J., Theory and
applications of fractional differential equations, North-Holland Mathematics
Studies, Vol.204, Elsevier, Amsterdam, 2006.

\item Samko, S.G., Kilbas, A.A. and Marichev, O., Farctional integrals and
derivaties, Gorden \&Breach, Berlin, 1993.

\item Kilbas, A.A., Trujillo, J.J., Differential equations of fractional
order: methods, results and problems I, Appl. Anal., 78, 153-192, 2001.

\item Oldham, K.B. and Spanier, J., The Fractional Calculus, Academic Press,
New York, 1974.

\item Podlubny, I., Fractional differential equations, Academic Press, New
York, 1999.

\item Miller, K.S. and Ross, B., An introduction to the fractional calculus
and fractional differntial Equations, Wiley, New York, 1993.

\item Hilfer, R., Applications of Fractional Calculus in Physics, World
Scientific, Singapore, 2000.

\item Levitan, B.M. and Sargsjan, I.S., Introduction to Spectral Theory:
Selfadjoint Ordinary Differential Operators, American Math. Soc., Pro.,
R.I., 1975.

\item Zettl, A., Sturm-Liouville Theory, Mathematical Surveys and
Monographs, American Mathematical Society, 2005.

\item Amrein, W.O., Hinz, A.M., Sturm-Liouville Theory; Past and Present,
Birkhauser, Basel, Switzerland, 2005.

\item Klimek, M. and Argawal, O. P., On a Regular Fractional Sturm-Liouville
Problem with Derivatives of Order in (0,1), 13th International Carpathian
Control Conference, Vysoke Tatry (Podbanske), Slovakia, May 28-31, 2012.

\item Klimek, M. and Argawal, O. P., Fractional Sturm-Liouville problem,
Computers and Mathematics with Applications, 66, 795-812, 2013.

\item Bas, E., Fundamental Spectral Theory of Fractional Singular
Sturm-Liouville Operator, Journal of Func. Spaces and Appl., Article ID
915830, 7 pages, 2013.

\item Bas, E., Metin, F., Fractional Singular Sturm-Liouville Operator for
Coulomb Potential, Advances in Differences Equations, DOI:
10.1186/1687-1847-2013-300, 2013.

\item Bas, E., Metin, F., Spectral Analysis for Fractional Hydrogen Atom
Equation, Advances in Pure Mathematics, 5, 767-773, 2015.

\item Ansari, A., On finite fractional Sturm-Liouville transforms, Integral
Transforms and Special Functions, 26(1) (2015), 51-64, 2015.

\item Ciesielski, M., Klimek, M. Blaszczyk, T., The fractional Sturm
Liouville problem Numerical approximation and application in fractional
diffusion, Journal of Computational and Applied Mathematics, 317, 573-588,
2017.

\item Zayernouri, M., Karniadakis, G. E., Fractional Sturm-Liouville
eigen-problems:Theory and numerical approximation, Journal of Computational
Physics, 252(1), 495-517, 2013.

\item Chai, G.,\textbf{\ }Positive solutions for boundary value problem of
fractional differential equation with p-Laplacian operator, \textit{Boundary
Value Problems, }18\textbf{,} 1-20, 2012.

\item Chen, T., Liu, W., Hu, Z., A boundary value problem for fractional
differential equation with p-Laplacian operator at resonance, \textit{%
Nonlinear Analysis, }75\textbf{, }3210-3217, 2012.

\item Chen, T., Liu, W., An anti-periodic boundary value problem for the
fractional differential equation with a p-Laplacian operator, \textit{%
Applied Mathematics Letters, }25\textbf{,} 1671-1675, 2012.

\item Mahmudov, N., Unul, S., Existence of solutions of fractional boundary
value problems with p-Laplacian operator, Boundary Value Problems, DOI
10.1186/s13661-015-0358-9.

\item Liu, X., Jia, M., Xiang, X., On the solvability of a fractional
differential equation model involving the p-Laplacian operator, Computers
and mathematics with Applications, 64, 3267-3275, 2012.

\item Lakshmikantham, V., Bainov, D.D., Simeonov, P.S., Theory of impulsive
differential equations, World Scientific, Singapore, 1989.

\item Samoilenko, A.M., Perestyuk, N.A., Impulsive differential equations,
World Scientific, Singapore, 1995.

\item Zavalishchin, S.T., Sesekin, A.N., Dynamic Impulse Systems: Theory and
Applications, Kluwer Academic Publishers Group, Dordrecht, 1997.

\item Zeidler, E., Nonlinear Functional Analysis and Its Applications-I:
Fixed-Point Theorems, Springer, New York, 1986.

\item Zhang, S., Positive solutions for baundary-value problems of nonlinear
fractional differential equations, Electronic Journal of Differential
Equations, 36, 1-12, 2006.

\item Tian, Y., Bai, Z., Existence results for the three-point impulsive
boundary value problem involving fractional differential equations,
Computers and Mathematics with Applications, 59, 2601-2609, 2010.

\item Abbas, S. and Benchohra, M., Impulsive partial hyperbolic functional
differential equations of fractional order with state-dependent delay,
Fractional Calculus \& Applied Analysis, 13, 225-244, 2010.

\item Agarwal, R.P., Benchohra, M. and Slimani, B.A., Existence results for
differential equations with fractional order and impulses,\ Memoirs on
Differential Equations and Mathematical Physics, 44, 1-21, 2008.

\item Ahmad, B. and Nieto, J.J., Existence of solutions for impulsive
anti-periodic boundary value problems of fractional order, Taiwanese Journal
of Mathematics, 15, 981-993, 2011.

\item Ahmad, B., Sivasundaram, S., Existence results for nonlinear impulsive
hybrid boundary value problems involving fractional differential equations,
Nonlinear Analysis: Hybrid Systems, 3, 251-258, 2009.

\item Zhou, J., Feng, M., Green's function for Sturm-Liouville-type boundary
value problems of fractional order impulsive differential equations and its
applications, Boundary Value Problems, 1, 69, 2014.
\end{enumerate}

\end{document}